\input amstex

\documentstyle{amsppt}
\loadbold

\def\<{\left<}										
\def\>{\right>}

\nologo

\def\qedd{
\hfill
\vrule height4pt width3pt depth2pt
\vskip .5cm
}

\magnification=\magstephalf

\topmatter
\title
The domain algebra of a $CP$-semigroup
\endtitle

\author William Arveson
\endauthor

\affil Department of Mathematics\\
University of California\\Berkeley CA 94720, USA
\endaffil

\date 18 May, 2000
\enddate
\thanks On appointment as a Miller Research 
Professor in the Miller Institute for Basic 
Research in Science.  
Support is also acknowledged from 
NSF grant DMS-9802474
\endthanks
%
%
%
\abstract 
A $CP$-semigroup (or {\it quantum dynamical semigroup}) 
is a semigroup $\phi=\{\phi_t: t\geq 0\}$ 
of normal completely positive linear maps on $\Cal B(H)$, 
$H$ being a separable Hilbert space, which satisfies 
$\phi_t(\bold 1)=\bold 1$ for all t and is 
continuous in the natural sense.

Let $\Cal D$ be the natural domain of the generator 
$L$ of $\phi$, $\phi_t=\exp{tL}$.  
Since the maps $\phi_t$ need not be multiplicative 
$\Cal D$ is typically an operator space, but not an algebra.  
However, we show that the set of operators 
$$
\Cal A=\{A\in\Cal D: A^*A\in\Cal D, AA^*\in\Cal D\}
$$
is a $*$-subalgebra of $\Cal B(H)$, indeed  $\Cal A$ is the 
largest self-adjoint algebra contained in $\Cal D$.  Because $\Cal A$ is 
a $*$-algebra one may consider its 
$*$-bimodule of noncommutative $2$-forms 
$\Omega^2(\Cal A)=\Omega^1(\Cal A)\otimes_\Cal A\Omega^1(\Cal A)$, 
and any linear mapping $L:\Cal A \to\Cal B(H)$ has a 
{\it symbol} $\sigma_L: \Omega^2(\Cal A)\to\Cal B(H)$, defined 
as a linear map by
$$
\sigma_L(a\,dx\,dy)=aL(xy)-axL(y)-aL(x)y+axL(\bold 1)y, \qquad a,x,y\in\Cal A.  
$$
The symbol is a homomorphism of $\Cal A$-bimodules for any 
$*$-algebra $\Cal A\subseteq\Cal B(H)$ and any linear map 
$L:\Cal A\to\Cal B(H)$.  
When $L$ is the generator of a $CP$-semigroup with domain
algebra $\Cal A$ above, we show that the symbol 
is negative in that 
$\sigma_L(\omega^*\omega)\leq 0$ for every $\omega\in\Omega^1(\Cal A)$
($-\sigma_L$ is in fact completely positive).  

Examples are given 
for which the domain algebra 
$\Cal A$ is, and is not, strongly dense in $\Cal B(H)$.  We also 
relate the generator of a $CP$-semigroup to its 
commutative paradigm, the Laplacian of a Riemannian manifold.  
\endabstract

\rightheadtext{the domain algebra}
\endtopmatter

\document

\subhead{1.  Basic properties of $\Cal A$}
\endsubhead
Let $\phi=\{\phi_t: t\geq 0\}$ be a $CP$-semigroup as 
defined in the abstract.  We first recall four characterizations 
of the domain of the generator of $\phi$.  
\proclaim{Lemma 1}
Let $A\in\Cal B(H)$.  The following are equivalent.  
\roster
\item"{(i)}" The limit 
$$
L(A)=\lim_{t\to 0+}\frac{1}{t}(\phi_t(A)-A)
$$
exists relative to the strong-$*$ topology of $\Cal B(H)$.  
\item"{(ii)}" the limit 
$$
L(A)=\lim_{t\to 0+}\frac{1}{t}(\phi_t(A)-A)
$$
exists relative to the weak operator topology of $\Cal B(H)$.  
\item"{(iii)}"
$$
\sup_{t>0}\frac{1}{t}\|\phi_t(A)-A\|\leq  M<\infty.  
$$
\item"{(iv)}" There is a sequence $t_n\to0+$ for which
$$
\sup_{n}\frac{1}{t_n}\|\phi_{t_n}(A)-A\|\leq  M<\infty.  
$$
\endroster
\endproclaim

\demo{proof}The implications (i) $\implies$ (ii) and 
(iii) $\implies$ (iv) are trivial, and 
(ii) $\implies$ (iii) is a straightforward consequence of the 
Banach-Steinhaus theorem.  

\demo{proof of (iv) $\implies$ (i)}
Since the unit ball of $\Cal B(H)$ is weakly sequentially compact, 
the hypothesis (iv) implies that 
there is a sequence $t_n\to0+$ such that 
$$
\frac{1}{t_n}(\phi_{t_n}(A)-A) \to T\in\Cal B(H)
$$
in the weak operator topology.  We claim: for every $s>0$, 
$$
\int_0^s \phi_\lambda(T)\,d\lambda = \phi_s(A) - A.  \tag{1.1}
$$
The integral on the left is interpreted as a weak integral; that 
is, for  $\xi,\eta\in H$, 
$$
\int_0^s\<\phi_\lambda(T)\xi,\eta\>\,d\lambda = 
\<\phi_s(A)\xi,\eta\>-\<A\xi,\eta\>.  
$$

To see that, fix $\lambda>0$.  Since 
$\phi_\lambda$ is weakly continuous on bounded 
sets in $\Cal B(H)$ we have 
$$
\frac{1}{t_n}(\phi_{\lambda+t_n}(A)-\phi_\lambda(A))=
\phi_\lambda(\frac{1}{t_n}(\phi_{t_n}(A)-A))\to \phi_\lambda(T)
$$
in the weak operator topology, as $n\to\infty$.  By the bounded convergence theorem, 
we find that for fixed $\xi,\eta\in H$, 
$$
\lim_{n\to\infty}\frac{1}{t_n}(\int_0^s\<\phi_{\lambda+t_n}(A)\xi,\eta\>\,d\lambda
-\int_0^s\<\phi_\lambda(A)\xi,\eta\>\,d\lambda) =
\int_0^s\<\phi_\lambda(T)\xi,\eta\>\,d\lambda.  
$$
Writing 
$$
\int_0^sf(\lambda+t_n)\,d\lambda-\int_0^sf(\lambda)\,d\lambda = 
\int_s^{s+t_n}f(\lambda)\,d\lambda - \int_0^{t_n}f(\lambda)\,d\lambda, 
$$
the left side of the preceding formula becomes 
$$
\lim_{n\to\infty}(\frac{1}{t_n}\int_s^{s+t_n}\<\phi_\lambda(A)\xi,\eta\>\,d\lambda -
\frac{1}{t_n}\int_0^{t_n}\<\phi_\lambda(A)\xi,\eta\>\,d\lambda)
$$
which, because of continuity of $\phi$ in the time parameter, is 
$\<\phi_s(A)\xi,\eta\>-\<A\xi,\eta\>$, as asserted in (1.1).  

To prove the strong-$*$ convergence asserted in (i), fix $\xi\in H$ and 
use (1.1) to write 
$$
\align
&\|\frac{1}{s}(\phi_s(A)\xi-A\xi)-T\xi\|= 
\frac{1}{s}\|\int_0^s\phi_\lambda(T)\xi\,d\lambda-\int_0^sT\xi\,d\lambda\| \\
&\leq \frac{1}{s}\int_0^s\|\phi_\lambda(T)\xi-T\xi\|\,d\lambda 
\leq (\frac{1}{s}\int_0^s\|\phi_\lambda(T)\xi-T\xi\|^2\,d\lambda)^{1/2}. 
\endalign
$$
The integrand of the last term expands as follows 
$$
\align
\|\phi_\lambda(T)\xi-T\xi\|^2&=\<\phi_\lambda(T)^*\phi_\lambda(T)\xi,\xi\> 
-2 \Re\<\phi_\lambda(T)\xi,T\xi\> +\|T\xi\|^2 \\
&\leq 
\<\phi_\lambda(T^*T)\xi,\xi\> 
-2 \Re\<\phi_\lambda(T)\xi,T\xi\> +\|T\xi\|^2, 
\endalign
$$
the last inequality by the Schwarz inequality for unital CP maps.  
Since $\phi_\lambda(T^*T)$ (resp. $\phi_\lambda(T)$)  tends weakly to 
$T^*T$ (resp. $T$) as $\lambda\to 0+$, it follows that 
$$
\limsup_{s\to 0+}\frac{1}{s}\int_0^s\|\phi_\lambda(T)\xi-T\xi\|^2\,d\lambda
\leq\<T^*T\xi,\xi\> 
-2 \<T\xi,T\xi\> +\|T\xi\|^2=0, 
$$
and we conclude that 
$\frac{1}{s}(\phi_s(A)-A)$ tends strongly to $T$ as $s\to0+$.  

Similarly,  
$\frac{1}{s}(\phi_s(A)-A)^*=\frac{1}{s}(\phi_s(A^*)-A^*)$ 
tends strongly to $T^*$.  \qedd
\enddemo

\proclaim{Definition}
Let $\Cal D$ be the set of all operators $A\in\Cal B(H)$ for which the 
four conditions of Lemma 1 are satisfied.  $L:\Cal D\to\Cal B(H)$ 
denotes the generator of $\phi$, 
$$
L(A) = \lim_{t\to0+}\frac{1}{t}(\phi_t(A)-A), \qquad A\in\Cal D.  
$$
\endproclaim
It is obvious that $\Cal D$ is a self-adjoint linear subspace of 
$\Cal B(H)$, that $L(A^*)=L(A)^*$ for $A\in\Cal D$, and a standard 
argument shows that $\Cal D$ is dense 
in $\Cal B(H)$ in the $\sigma$-strong 
operator topology.  

\proclaim{Lemma 2} For every operator $A\in\Cal D$ we have 
$$
\|L(A)\| = \sup_{t>0}\frac{1}{t}\|\phi_t(A)-A\|.   
$$
\endproclaim
\demo{proof} The inequality $\leq$ is clear from the fact that 
$L(A)$ is the weak limit of operators $\frac{1}{t}(\phi_t(A)-A)$ 
near $t=0+$, i.e., 
$$
\|L(A)\|\leq \limsup_{t\to0+}\frac{1}{t}\|\phi_t(A)-A\|\leq 
\sup_{t>0}\frac{1}{t}\|\phi_t(A)-A\|.  
$$

For $\geq$, set $T=L(A)$.  Using (1.1), we can write for every $t>0$
$$
\frac{1}{t}\|\phi_t(A)-A\|=
\frac{1}{t}\|\int_0^t\phi_\lambda(T)\,d\lambda\|\leq 
\frac{1}{t}\int_0^t\|\phi_\lambda(T)\|\,d\lambda\leq \|T\|,  
$$
since $\|\phi_\lambda\|\leq 1$ for every $\lambda\geq 0$.\qedd
\enddemo

\proclaim{Theorem A}
$\Cal A=\{A\in\Cal D: A^*A\in\Cal D, AA^*\in\Cal D\}$ is a 
$*$-subalgebra of $\Cal B(H)$.  
\endproclaim

\demo{proof}
$\Cal A$ is obviously a self-adjoint set of operators.  
We have to show that $\Cal A$ is a vector space satisfying 
$\Cal A\cdot\Cal A\subseteq \Cal A$.   

  Fix $t>0$.  By Stinespring's 
theorem we can write 
$$
\phi_t(X) = V_t^*\pi_t(X)V_t, \qquad X\in\Cal B(H) \tag{1.2}
$$
where 
$V_t$ is an isometry from $H$ into some other Hilbert space $H_t$ and 
where $\pi_t: \Cal B(H)\to \Cal B(H_t)$ is a {\it normal} $*$-homomorphism 
of von Neumann algebras.  $P_t = V_tV_t^*$ is a self-ajoint projection 
in $\Cal B(H_t)$.  

For $t>0$ we will consider the seminorms $p_t$, $q_t$ defined on 
$\Cal B(H)$ as follows
$$
\align
p_t(X) &= t^{-1}\|\phi_t(X)-X\|, \\
q_t(X) &= t^{-1/2}\|P_t\pi_t(X) - \pi_t(X)P_t\|, \qquad X\in\Cal B(H).  
\endalign
$$

\proclaim{Lemma 3} For every operator $X\in\Cal B(H)$ we have 
the following characterizations.  
\roster
\item"{(i)}" 
$X\in\Cal D$ iff 
$$
\sup_{t>0}p_t(X)<\infty,
$$
and in that case $\|L(X)\|=\sup_{t>0}p_t(X)$.  

\item"{(ii)}"  $X\in\Cal A$ iff both $\sup_{t>0}p_t(X)$ and 
$\sup_{t>0}q_t(X)$ are finite, and in that case 
$$
\max(\|\sigma_L(dX^*\,dX)\|^{1/2},\|\sigma_L(dX\,dX^*)\|^{1/2})
\leq \limsup_{t\to 0+}q_t(X),
$$
where $\sigma_L(dX^*\,dX)$ and $\sigma_L(dX\,dX^*)$ are the 
operators in $\Cal B(H)$ defined by
$$
\align
\sigma_L(dX^*\,dX)&=L(X^*X)-X^*L(X)-L(X^*)X, \\
\sigma_L(dX\,dX^*)&=L(XX^*)-XL(X^*)-L(X)X^*
\endalign
$$
\endroster
\endproclaim

\remark{Remark}
The second assertion of Lemma 3 requires clarification.  By definition, 
an operator $X$ belongs to $\Cal A$ iff all four operators 
$X,X^*, X^*X, XX^*$ belong to the domain of the generator $L$ of 
$\phi=\{\phi_t: t\geq 0\}$.  In that case both operators 
$\sigma_L(dX^*\,dX)$ and $\sigma_L(dX\,dX^*)$
are well defined by the above formulas. The ``symbol" map $\sigma_L$ 
will be discussed more fully in section 3.
\endremark

\demo{proof of Lemma 3}  The assertion (i) follows from Lemmas 
1 and 2 above.  In order to prove (ii) we require
the following more concrete expression for the seminorm $q_t$,
$$
\align
&q_t(X) =\\
&\max(\|\frac{1}{t}(\phi_t(X^*X)-\phi_t(X)^*\phi_t(X))\|^{1/2}, 
\|\frac{1}{t}(\phi_t(XX^*)-\phi_t(X)^*\phi_t(X^*))\|^{1/2}).  \tag{1.3}
\endalign
$$ 

To prove (1.3) we decompose the commutator $\pi_t(X)P_t-P_t\pi_t(X)$ 
into a sum
$$
\pi_t(X)P_t-P_t\pi_t(X)=(\bold 1-P_t)\pi_t(X)P_t - P_t\pi_t(X)(\bold 1-P_t).  
$$
Since the first term $(\bold 1-P_t)\pi_t(X)P_t$ has 
initial space in $P_tH_t$ and final space in $(\bold 1-P_t)$, and the 
second term has the opposite property, it follows that 
$$
\|\pi_t(X)P_t-P_t\pi_t(X)\|=
\max(\|(\bold 1-P_t)\pi_t(X)P_t\|, \|P_t\pi_t(X)(\bold 1-P_t)\|).  
$$
We have
$$
\align
\|(\bold 1-P_t)\pi_t(X)P_t\|^2 &=
\|V_t^*\pi_t(X^*)(\bold 1-P_t)\pi_t(X)V_t\| \\
&=\|V_t^*\pi_t(X^*X)V_t-V_t^*\pi_t(X^*)V_tV_t^*\pi_t(X)V_t\| \\
&=\|\phi_t(X^*X)-\phi_t(X)^*\phi_t(X)\|. 
\endalign 
$$
Similarly,
$$
\|P_t\pi_t(X)(\bold 1-P_t)\|^2 =\|V_t^*\pi_t(X)(\bold 1-P_t)\pi_t(X^*)V_t\|=
\|\phi_t(XX^*)-\phi_t(X)^*\phi_t(X^*)\|, 
$$
and formula (1.3) follows from these two expressions.  

Now if $X\in\Cal A$ then  all four operators $X, X^*,X^*X, XX^*$
belong to $\Cal D$, hence all four limits 
$$
\align
&\lim_{t\to0+}\frac{1}{t}(\phi_t(X^*X)-X^*X)=L(X^*X), \\
&\lim_{t\to0+}\frac{1}{t}(\phi_t(XX^*)-X^*X)=L(XX^*), \\
&\lim_{t\to0+}\frac{1}{t}(\phi_t(X)-X)=L(X), \\
&\lim_{t\to0+}\frac{1}{t}(\phi_t(X^*)-X^*)=L(X^*)
\endalign
$$
exist relative to the strong operator topology.  Writing 
$$
\align
&\phi_t(X^*X)-\phi_t(X)^*\phi_t(X)= \\
(\phi_t(X^*X)-&X^*X) - X^*(\phi_t(X)-X) -(\phi_t(X^*)-X^*)\phi_t(X)\tag{1.4}
\endalign
$$
and using strong continuity of multiplication on bounded sets, 
we find that the limit 
$$
\lim_{t\to0+}\frac{1}{t}(\phi_t(X^*X)-\phi_t(X^*)\phi_t(X)) = 
L(X^*X)-X^*L(X)-L(X^*)X = \sigma_L(dX^*\,dX)
$$
exists relative to the strong operator topology.  

In the same way we deduce the existence of the strong limit
$$
\lim_{t\to0+}\frac{1}{t}(\phi_t(XX^*)-\phi_t(X)\phi_t(X^*)) = 
L(XX^*)-XL(X^*)-L(X)X^* = \sigma_L(dX\,dX^*).
$$
It follows that for every $X\in \Cal A$ 
the seminorms $q_t(X)$ are bounded for $t>0$, 
and for such $X$ we have 
$$
\max(\|\sigma_L(dX^*\,dX)\|^{1/2},\|
\sigma_L(dX\,dX^*)\|^{1/2})\leq \limsup_{t\to0+}q_t(X).  
$$

Conversely, suppose we are given an operator $X\in\Cal D$ for 
which the seminorms $q_t(X)$ are bounded for $t>0$.  We have to show 
that $X^*X$ and $XX^*$ belong to $\Cal D$; since $\Cal D$ 
is self-adjoint and the seminorms $q_t$ are symmetric 
in that $q_t(X^*)=q_t(X)$, it is enough to show that 
$X^*X$ belong to $\Cal D$.  
(1.4) implies that for fixed $t>0$,
$$
\align
&\phi_t(X^*X)-X^*X= \\
(\phi_t(X^*X)-&\phi_t(X^*)\phi_t(X)) + X^*(\phi_t(X)-X) +
(\phi_t(X^*)-X^*)\phi_t(X)\tag{1.5}
\endalign
$$
Because of (1.3), the first term on the right of (1.5) is bounded 
in norm by $M_1\cdot t$ where $M_1$ is a positive constant.  
Similarly, since  $X$ and 
$X^*$ belong to $\Cal D$ the second and third terms are bounded in norm 
by terms of the form $M_2\cdot t$ and $M_3\cdot t$ respectively, hence 
$$
\|\phi_t(X^*X)-X^*X\|\leq (M_1+M_2+M_3)\cdot t.  
$$
By Lemma 1, $X^*X$ must belong to $\Cal D$.  
\qedd
\enddemo

Turning now to the proof of Theorem A, (or more properly, to the proof 
that $\Cal A$ is an algebra), Lemma 3 tells us that $\Cal A$ consists 
of all operators $X\in\Cal B(H)$ for which 
$$
\sup_{t>0}p_t(X)<\infty, \qquad {\text{and }}\sup_{t>0}q_t(X)<\infty.  
$$
Since $p_t$ and $q_t$ are both seminorms, it follows that 
$\Cal A$ is a complex vector space which is obviously closed 
under the $*$-operation.  

To see that $\Cal A$ is closed under multiplication, 
pick  $X,Y\in\Cal A$.  According to Lemma 3, it is enough to show
$$
\sup_{t>0}q_t(XY)<\infty \tag{1.6}
$$
and 
$$
\sup_{t>0}p_t(XY)<\infty \tag{1.7}
$$
To prove (1.6) we claim that 
$$
q_t(XY) \leq q_t(X)\|Y\| + \|X\| q_t(Y).  \tag{1.8}
$$
Indeed, writing $[A,B]$ for the commutator $AB-BA$ we have 
$$
[P_t, \pi_t(XY)]=[P_t,\pi_t(X)]\pi_t(Y) + \pi_t(X)[P_t,\pi_t(Y)], 
$$
and hence 
$$
\align
q_t(XY)&=t^{-1/2}\|[P_t,\pi_t(XY)]\|  \\
&\leq 
t^{-1/2}\|[P_t,\pi_t(X)\|\cdot\|\pi_t(Y)\| + 
\|\pi_t(X)\|\cdot t^{-1/2}\|[P_t,\pi_t(Y)]\|, 
\endalign
$$
from which (1.8) is evident.  

Finally, consider the condition (1.7).  By definition of $\Cal A$,
$A\in\Cal A$ implies $A^*A\in\Cal D$.  Since $\Cal A$ is 
now known to be a linear space we can assert 
that if $X,Y\in\Cal A$ then for every 
$k=0,1,2,3$ we have $Y+i^kX\in\Cal A$, 
hence $(Y+i^kX)^*(Y+i^kX)\in\Cal D$ and 
by the polarization formula 
$$
X^*Y=\frac{1}{4}\sum_{k=0}^3 i^k(Y+i^kX)^*(Y+i^kX), 
$$
$X^*Y$ must also belong to $\Cal D$.  
Since $\Cal A^*=\Cal A$, we can 
replace $X^*$ with $X$ to conclude that $XY\in\Cal D$.  
(1.7) now follows from Lemma 3 (i).  \qedd
\enddemo
\enddemo

\proclaim{Corollary}
Let $\Cal D$ be the domain of the generator of a $CP$-semigroup acting
on $\Cal B(H)$ and 
let $A$ be a self-adjoint operator such that $A\in\Cal D$ and $A^2\in\Cal D$.  
Then $p(A)\in\Cal D$ for every polynomial $p(x)=a_0+a_1x+\dots+a_nx^n$.  
\endproclaim

\subhead{2.  Examples}
\endsubhead
In this section we describe two classes of examples which 
are in a sense at opposite extremes.  In the first class of 
examples of $CP$-semigroups $\phi=\{\phi_t: t\geq 0\}$, 
each $\phi_t$ leaves the $C^*$-algebra $\Cal K$ of all 
compact operators invariant,
$\phi_t(\Cal K)\subseteq \Cal K$, its domain algebra 
$\Cal A$ is strongly dense in $\Cal B(H)$, and 
its generator restricts to a ``second 
order" differential operator on $\Cal A$ (see formula 
(1.1) of \cite{1}).  In the second class of examples, the 
individual maps satisfy $\phi_t(\Cal K)\cap\Cal K=\{0\}$ for $t>0$, 
$\Cal A$ is not strongly dense in $\Cal B(H)$, 
and its generator is degenerate in the sense that it 
restricts to a {\it derivation} on $\Cal A$.  

We first recall the class of examples of $CP$-semigroups 
of \cite{1}, including the heat flow of the $CCR$ algebra. 
While for simplicity we confine the discussion to 
the case of one degree of freedom, the 
reader will note that everything carries over 
verbatim to the case of $n$ degrees of freedom, $n=1,2,\dots$.  
 
Let $\{W_z: z\in\Bbb R^2\}$ be an irreducible Weyl system 
acting on a Hilbert space $H$.  Thus, $z\in\Bbb R^2\mapsto W_z$
is a strongly continuous mapping from $\Bbb R^2$ into the unitary 
operators on $H$ which satisfies the canonical commutation 
relations in Weyl's form 
$$
W_{z_1}W_{z_2}=e^{i\omega(z_1,z_2)}W_{z_1+z_2}, 
\qquad z_1,z_2\in\Bbb R^2,
$$
$\omega$ denoting the symplectic form on $\Bbb R^2$ given 
by 
$$
\omega((x,y), (x^\prime,y^\prime))= 
\frac{1}{2}(x^\prime y-xy^\prime).  
$$
Let $\{\mu_t: t\geq 0\}$ be a one-parameter family of probability 
measures on $\Bbb R^2$ which is a semigroup under the 
natural convolution of measures
$$
\mu * \nu(S)=\int_{\Bbb R^2\times\Bbb R^2}\chi_S(z+w)\,d\mu(z)\,d\nu(w),
$$ 
which satisfies $\mu_0=\delta_{(0,0)}$, 
and which is measurable in $t$ in the natural sense.  
It is convenient to define 
the Fourier transform of a measure $\mu$ in terms of the symplectic 
form $\omega$ as follows, 
$$
\hat\mu(z)=\int_{\Bbb R^2} e^{i\omega(z,\zeta)}\,d\mu(\zeta), 
\qquad z\in \Bbb R^2.  
$$
Given such a semigroup of probability measures $\{\mu_t: t\geq 0\}$
there is a unique $CP$ semigroup $\phi=\{\phi_t: t\geq 0\}$ acting 
on $\Cal B(H)$ which satisfies 
$$
\phi_t(W_z)=\hat\mu_t(z)W_z, \qquad z\in \Bbb R^2, \quad t\geq 0
$$
see \cite{1}, Proposition 1.7.  
Two cases of particular interest are

$$
\phi_t(W_z)=e^{-t|z|^2}W_z, \qquad t\geq 0\tag{CCR heat flow}
$$
where $|(x,y)|$ denotes the Euclidean norm $(x^2+y^2)^{1/2}$, and 
$$
\phi_t(W_z)=e^{-t|z|}W_z,\qquad t\geq 0. \tag{Cauchy flow}
$$
For both of these examples a straightforward 
estimate shows that for 
fixed $z\in\Bbb R^2$ there is a constant $M>0$ such that 
$$
\|\phi_t(W_z)-W_z\|=|\hat\mu_t(z)-1|\leq M\cdot t, \qquad t>0
$$
and hence $W_z\in\Cal D$.  Since $W_z$ is unitary, 
$\bold 1=W_z^*W_z=W_zW_z^*$ belongs to $\Cal D$ , 
and hence $W_z$ belongs to the domain algebra 
$\Cal A$ of $\phi$ for every $z\in\Bbb R^2$.  We conclude that 
for these examples, the domain algebra is strongly dense 
in $\Cal B(H)$.  

Indeed, it is not hard to show that $\Cal A$ contains a 
$*$-algebra of compact operators that is norm-dense in 
the algebra $\Cal K$ of all compact operators.  Unlike 
the examples to follow, these flows 
leave $\Cal K$ invariant in the sense that 
$\phi_t(\Cal K)\subseteq \Cal K$ for all $t\geq 0$, and 
can therefore be considered as $CP$-semigroups which act 
on the separable $C^*$-algebra $\Cal K$,  rather than 
than as $CP$-semigroups acting on $\Cal B(H)$.  

\vskip0.1truein
We now describe a class of examples 
of $CP$ semigroups whose 
domain algebras are {\it not} strongly 
dense in $\Cal B(H)$.  
These examples are inspired by a class of $CP$ semigroups
that have emerged in recent work of Robert Powers, to whom 
we are indebted for useful discussions.  

Let $H=L^2(0,\infty)$ and let $U=\{U_t: t\geq 0\}$ be 
the semigroup of isometries $U_t \xi(x)=\xi(x-t)$ for 
$x\geq t$, $U_t\xi(x)=0$ for $0\leq x<t$.  Fix a 
real number $\alpha>0$, and let $f$ be the unit 
vector in $L^2(0,\infty)$ obtained by normalizing the exponential 
function $u(x)= e^{-\alpha x}$, $x\geq0$.  One has 
$U_t^*f=e^{-\alpha t}f$ for every $t\geq 0$, hence 
the vector state $\omega(A)=\<Af,f\>$ satisfies 
$\omega(U_tAU_t^*)=e^{-2\alpha t}\omega(A)$, 
$A\in\Cal B(H)$.  

We consider the family of unit-preserving normal completely positive 
maps $\phi=\{\phi_t: t\geq0\}$ defined on $\Cal B(H)$ by 

$$
\phi_t(A) = \omega(A)E_t + U_tAU_t^*, \qquad t\geq 0.  
$$
where $E_t=\bold 1-U_tU_t^*$ is the projection on 
the subspace $L^2(0,t)\subseteq L^2(0,\infty)$.  Since 
$$
\omega(E_t) = \omega(\bold 1)-\omega(U_tU_t^*)=1-e^{-2\alpha t}, 
$$
it follows that $\omega(\phi_t(A))=\omega(A)$ for 
every $A$.  
A routine computation now shows that $\phi$ satisfies the 
semigroup property $\phi_s\circ\phi_t=\phi_{s+t}$, hence 
$\phi$ is a $CP$ semigroup.

Let $\Cal D$ be the domain of the generator of $\phi$ and 
let $\Cal A$ be the domain algebra
$$
\Cal A=\{A\in\Cal D: A^*A\in\Cal D, AA^*\in\Cal D\}.  
$$
Theorem A implies that $\Cal A$ is a unital $*$-algebra
and we calculate its strong closure.  

\proclaim{Proposition}
The strong closure of $\Cal A$ consists of all operators 
$B\in\Cal B(H)$ such that $B$ commutes with the rank-one 
projection $f\otimes\bar f$.  
\endproclaim

Thus the strong closure of 
$\Cal A$ consists of all operators $B$ 
such that both $B$ and $B^*$ have $f$ as an eigenvector.  

\demo{proof} By Lemma 1, the domain $\Cal D$ of the 
generator of $\phi$ consists of all operators $A$ 
with the property 
$$
\|\phi_t(A)-A\|\leq M\cdot t, \qquad {\text{for all }}t\geq 0, \tag{2.1}
$$
where $M$ is a positive constant depending on $A$.  

First, we show that $f\otimes f$ commutes with $\Cal A$.  
Choose $A\in\Cal A$.  In order to show that $A$ commutes 
with $f\otimes \bar f$, it is enough to show that 
$$
\omega(A^*A)=\omega(AA^*)=|\omega(A)|^2,\tag{2.2}
$$ since (2.2) implies 
$$
\|Af-\omega(A)f\|^2=\omega(A^*A)-2|\omega(A)|^2+|\omega(A)|^2=0,
$$
and similarly $\|A^*f-\omega(A^*)f\|=0$.  
Multiplying $\phi_t(A)-A$ on the right by $E_t$ and using the fact 
that $\phi_t(A)E_t=\omega(A)E_t$ we conclude from (2.1) that 
$$
\lim_{t\to0}\|\omega(A)E_t-AE_t\|=0.  
$$
Replacing $A$ with $A^*A$ and $AA^*$ one also finds 
$$
\lim_{t\to0}\|\omega(A^*A)E_t-A^*AE_t\|=
\lim_{t\to0}\|\omega(AA^*)E_t-AA^*E_t\|=0.  
$$
Taken together, these three limits imply that 
$\omega(A^*A)=\omega(AA^*)=|\omega(A)|^2$, as required. 

To prove the opposite inclusion
it is enough to show that for every self-adjoint operator 
$A\in\Cal B(H)$ satisfying $Af=0$ there is a sequence $A_n$ 
of self-adjoint operators in $\Cal A$ which converges weakly  
to $A$ (recall that $\Cal A$ is a self-adjoint algebra 
containing the identity).  Fix such an $A$ and, for every $\epsilon>0$, set 
$$
A_\epsilon = \phi_\epsilon(A)=
\omega(A)E_\epsilon+U_\epsilon AU_\epsilon^*=
U_\epsilon AU_\epsilon^*.  
$$
$A_\epsilon$ converges weakly to $A$ as $\epsilon\to 0$.  
Moreover, $A_\epsilon$ is supported in the 
interval $(\epsilon,\infty)$ in 
the sense that $A_\epsilon E_\epsilon=E_\epsilon A_\epsilon=0$, 
and in addition we have $A_\epsilon f=0$ since 
$$
A_\epsilon f=U_\epsilon AU_\epsilon^* f=
e^{-\alpha\epsilon}U_\epsilon Af = 0.  
$$
We show that each $A_\epsilon$ can be 
weakly approximated by self-adjoint elements of 
the domain algebra.  

\proclaim{Lemma} Suppose $\epsilon>0$ and let 
$A$ be a self-adjoint operator in $\Cal B(H)$ such that (i) 
$Af=0$ and (ii) $A$ is supported in $(\epsilon,\infty)$ 
in the sense that $AE_\epsilon=E_\epsilon A=0$.  
Let $u$ be a $C^\infty$ 
function having compact support in $[0,\epsilon]$ and consider 
$$
B = \int_0^\infty u(s)U_sAU_s^*\,ds=
\int_0^\infty u(s)\phi_s(A)\,ds.  
$$

Then $B^n\in\Cal D$ for every $n=1,2,\dots$, and in particular 
$B\in\Cal A$.  
\endproclaim

\demo{proof}
Observe first that $B$ has both properties (i) and (ii), 
hence so does $B^n$ for every $n$.  
Thus for $t<\epsilon$ we have 
$$
\phi_t(B^n)-B^n=U_tB^nU_t^*-B^n=U_tB^nU_t^*-B^nU_tU_t^*=
(U_tB^n-B^nU_t)U_t^*.  
$$
This implies that for sufficiently small $t$
$$
\|\phi_t(B^n)-B^n\|=\|U_tB^n-B^nU_t\|.  
$$
We conclude that $B^n\in\Cal D$ iff there is a constant 
$K>0$ such that 
$$
\|U_tB^n-B^nU_t\|\leq K\cdot t, \qquad {\text{for all }}t>0.  \tag{2.3}
$$
To prove (2.3), one uses the Leibniz rule for the derivation 
$D(X) = U_tX-XU_t$ to estimate $\|U_tB^n-B^nU_t\|$ in 
terms of $\|U_tB-BU_t\|$, 
$$
\|D(B^n)\|\leq n\cdot\|B\|^{n-1}\|D(B)\|=n\cdot\|B\|^{n-1}\|U_tB-BU_t\|.  
$$
Since $B$ has been smoothed it belongs to the domain 
$\Cal D$, hence there is a constant $M$ such that 
$\|U_tB-BU_t\|\leq M\cdot t$, hence 
$\|U_tB^n-B^nU_t\|\leq nM\|B\|^{n-1}\cdot t$.  \qedd
\enddemo

The proof of the Proposition is completed by choosing 
$A=A_\epsilon$ in the hypothesis of the Lemma and by 
choosing a sequence $u_k$ 
of nonnegative $C^\infty$ functions, each of which has 
integral $1$, such that $u_k(x)=0$ outside the interval 
$0\leq x\leq 1/k$.  A standard argument shows that the 
sequence of self-adjoint operators 
$$
B_k=\int_0^\infty u_k(s)\phi_s(A_\epsilon)\,ds
$$
converges weakly to $A_\epsilon$, and the Lemma implies that 
$B_k\in\Cal A$ for $k> 1/\epsilon$.
\qedd
\enddemo

Thus the strong closure $\Cal A^-$ of $\Cal A$ has the form 
$\Cal B(H_0)\oplus\Bbb C$ where $H_0\subseteq H$ is 
a subspace of codimension one in $H$, and 
the following implies that 
these examples are 
``almost" $E_0$-semigroups in the sense that there 
is an $E_0$-semigroup $\alpha=\{\alpha_t: t\geq 0\}$ 
acting on $\Cal B(H_0)$ such that $\phi_t$ 
acts as follows on $\Cal A^-$,
$$
\phi_t(B\oplus\lambda)=\alpha_t(B)\oplus\lambda, 
\qquad B\in\Cal B(H_0),\quad\lambda\in\Bbb C.    
$$

\proclaim{Corollary}
Let $\bar\Cal A$ be the strong closure of $\Cal A$.  Then 
$\phi_t(\bar\Cal A)\subseteq \bar\Cal A$ for every $t\geq0$ 
and $\{\phi_t\restriction_{\bar\Cal A}: t\geq 0\}$ 
is a semigroup of endomorphisms of this von Neumann algebra.  
\endproclaim

\demo{proof}
We show that 
$\phi_t(\Cal A)\subseteq\bar\Cal A$, and for $A,B\in\Cal A$ 
one has $\phi_t(AB)=\phi_t(A)\phi_t(B)$.  

Choose $A\in\Cal A$, and let $f$ and $\omega(A)=\<Af,f\>$ 
be as in the definition of $\phi_t$,
$$
\phi_t(A)=\omega(A)E_t+U_tAU_t^*, \qquad A\in\Cal A, 
\quad t\geq 0.  
$$
Since $f$ is an eigenvector for both $A$ and $A^*$ 
and $U_t^*f=e^{-\alpha t}f$, one can verify 
directly that $\phi_t(A)f=\omega(A)f$ and $\phi_t(A)^*f=\omega(A^*)f$, 
and the Proposition implies that $\phi_t(\Cal A)\subseteq\bar\Cal A$.  
Finally, for $A$, $B\in\Cal A$ one has $\omega(AB)=\omega(A)\omega(B)$, and 
$\phi_t(AB)=\omega(AB)E_t+U_tABU_t^*=\omega(A)\omega(B)E_t+U_tAU_t^*U_tBU_t^*
=\phi_t(A)\phi_t(B)$.  By normality of $\phi_t$, the 
formula $\phi_t(AB)=\phi_t(A)\phi_t(B)$ persists for operators $A$,$B$ in 
the strong closure of $\Cal A$.\qedd

\enddemo

\subhead{3.  The symbol of the generator: properties and structure}
\endsubhead

There are two useful characterizations of the generators 
of {\it uniformly continuous} $CP$-semigroups, i.e., those 
whose generators are everywhere defined bounded linear 
maps on $\Cal B(H)$.  The first is due to Lindblad \cite{24} 
and independently to Gorini {\it et al} \cite{20}
(also see \cite{13}, Theorem 4.2).  The second characterization
is due to Evans and Lewis \cite{19},
based on work of Evans \cite{16}.  These two results can 
be paraphrased as follows.  

\proclaim{Theorem}
Let $L: \Cal B(H)\to\Cal B(H)$ be a bounded linear map and 
let $\phi=\{\phi_t:t\geq0\}$ be the semigroup defined 
on $\Cal B(H)$ by $\phi_t=\exp(tL)$.  The following are 
equivalent.  
\roster
\item $\phi_t$ is a completely positive map for every $t\geq 0$.  
\item (Lindblad, Gorini et al)
$L$ admits a decomposition 
$$
L(A)=P(A) + BA+AB^*, \qquad A\in\Cal B(H)
$$
where $P$ is a completely positive linear map and $B\in\Cal B(H)$.  
\item
(Evans and Lewis)  For every finite set of operators 
$A_1,\dots, A_n, B_1,\dots, B_n\in\Cal B(H)$ which satisfy 
$A_1B_1+\dots+A_nB_n=0$, we have 
$$
\sum_{i,j=1}^nA_j^*L(B_j^*B_i)A_i\geq 0.  
$$
\endroster
\endproclaim

A linear map $L:\Cal B(H)\to\Cal B(H)$ satisfying property 
(3) of Theorem 3.1 is called {\it conditionally completely 
positive} \cite{17}.  While the characterization (2) tells 
us exactly which bounded linear maps generate $CP$ semigroups, 
the cited decomposition of $L$ into a sum of more 
familiar mappings is unfortunately not unique.  

The purpose of this section is to make two observations.  
First, we point out that the notion of a conditionally completely 
positive linear map defined on a $*$-algebra is more 
properly formulated in terms of the bimodule of 
noncommutative 2-forms over that algebra; and once that 
is done the ``symbol" of the map becomes analogous 
to a Riemannian metric.  Second, we 
show that by making use of the domain algebra of section 1, 
this notion becomes appropriate for the 
generators of {\it arbitrary} $CP$-semigroups.

Let $\Cal A$ be the domain algebra of a $CP$ semigroup
$\phi=\{\phi_t: t\geq 0\}$ acting on $\Cal B(H)$
$$
\Cal A = \{A\in\Cal D: A^*A\in\Cal D, AA^*\in\Cal D\}, 
$$
where $\Cal D$ is the natural domain of the generator 
$L$ of $\phi$.  We first recall the definition of the module of 
noncommutative 1-forms 
$\Omega^1(\Cal A)$, and 2-forms $\Omega^2(\Cal A)$.  The algebraic 
tensor product of vector spaces $\Cal A\otimes\Cal A$ can 
be considered an involutive bimodule over $\Cal A$, with 
$$
\align
a(x\otimes y)b&=ax\otimes yb,\\
(x\otimes y)^*&=y^*\otimes x^*.  
\endalign
$$
The map $d: \Cal A\to\Cal A\otimes \Cal A$ defined by 
$dx=\bold 1\otimes x-x\otimes\bold 1$ is a derivation for 
which $(dx)^*=-d(x^*)$, and it is a universal derivation of 
$\Cal A$ in the sense that if $E$ is any $\Cal A$-bimodule 
and $D:\Cal A\to E$ is a linear map satisfying 
$D(xy)=xD(y)+D(x)y$ for all $x,y\in\Cal A$, then there is a 
unique homomorphism of $\Cal A$-modules 
$\theta: \Omega^1(\Cal A)\to E$ such that $\theta\circ d=D$.  
Every element of $\Omega^1(\Cal A)$ is a finite sum of the 
form
$$
\omega=\sum_{k=1}^ra_k\,dx_k, 
$$
and the involution in $\Omega^1(\Cal A)$ satisfies 
$$
(a\,dx)^*=-d(x^*)a^*=-d(x^*a^*)+x^*\,d(a^*).  
$$
Finally, $\Omega^1(\Cal A)$ is the kernel of the multiplication 
map $\mu: \Cal A\otimes\Cal A\to \Cal A$ defined by 
$\mu(x\otimes y)=xy$, and thus we have an exact sequence 
of $\Cal A$-modules 
$$
0\longrightarrow \Omega^1(\Cal A)\longrightarrow 
\Cal A\otimes \Cal A\underset\mu\to\longrightarrow\Cal A\longrightarrow 0.  
\tag{3.1}
$$

$\Omega^2(\Cal A)$ is defined by 
$$
\Omega^2(\Cal A)=\Omega^1(\Cal A)\otimes_\Cal A \Omega^1(\Cal A), 
$$
and any element of $\Omega^2(\Cal A)$ can be written as a sum 
$$
\omega=\sum_{k=1}^r a_k\,dx_k\,dy_k.  
$$
The involution in $\Omega^2(\Cal A)$ satisfies 
$$
(a\,dx\,dy)^*=d(y^*)\,d(x^*)\,a^*=
d(y^*)\,d(x^*a^*)-d(y^*x^*)d(a^*)+y^*d(x^*)d(a^*).  
$$

Since $\Cal A$ is a $*$-subalgebra of $\Cal B(H)$, we may also 
think of $\Cal B(H)$ as an $\Cal A$-bimodule.  Now a straightforward 
argument shows that for every linear mapping 
$L: \Cal A\to\Cal B(H)$ there is a unique homomorphism of 
bimodules $\sigma_L: \Omega^2(\Cal A)\to\Cal B(H)$ which satisfies 
$$
\sigma_L(dx\,dy)=L(xy)-xL(y)-L(x)y+xL(\bold 1)y, 
\qquad x,y\in \Cal A  \tag{3.2}
$$
(see section 2 of \cite{4} for more detail).  
$\sigma_L\in\hom(\Omega^2,\Cal B(H))$ is called the {\bf symbol} 
of the linear map $L$.  

Consider now the special case in which $\Cal A=\Cal B(H)$,  
 $L:\Cal B(H)\to\Cal B(H)$ is a {\it bounded} linear 
mapping, and let $\phi=\{\phi_t=\exp{tL}: t\geq 0\}$ is the semigroup
of bounded operators on $\Cal B(H)$ generated by $L$.  The 
preceding theorem gives two characterizations of the  maps 
$L$ for which each $\phi_t=\exp{tL}$ is completely positive; 
however, the following characterization is perhaps more in spirit with 
the theory of differential operators on manifolds.  

\proclaim{Theorem}
Let $L:\Cal B(H)\to\Cal B(H)$ be a bounded linear map.  To the two 
characterizations (2) (3) above, one can append the following equivalent 
condition
\roster
\item"{(4)}"  The symbol $\sigma_L: \Omega^2(\Cal B(H))\to\Cal B(H)$ 
satisfies 
$$
\sigma_L(\omega^*\omega)\leq 0, \qquad 
{\text{for every }}\omega\in\Omega^1(\Cal B(H)).
$$
\endroster
\endproclaim

This characterization is Proposition 1.6 of \cite{5}; a fuller 
discussion of these issues can be found in \cite{4}.  Notice that 
the sense of the inequality $\leq$ is determined by the fact that 
the involution in $\Omega^1$ satisfies $(dx)^*=-d(x^*)$, and hence 
for $\omega=dx$ we have $\omega^*\omega=-d(x^*)\,dx$.  In particular, 
for $\omega=dx$ where $x$ is a self-adjoint element we have 
$\sigma_L(\omega^2)\geq 0$ while $\sigma_L(\omega^*\omega)\leq 0$.  

\remark{Remarks}
There is a rather compelling analogy between this characterization 
of the generators of CP semigroups and the generator of the 
heat flow of a Riemannian manifold, namely the Laplacian.  
More precisely, let $M$ be a complete (but not necessarily 
compact) Riemannian manifold and consider 
its natural Hilbert space $L^2(M)$.  The Laplacian $\Delta$ acts naturally 
as a densely defined operator on $L^2(M)$ and 
generates a semigroup of bounded operators $\exp{t\Delta}$, $t\geq 0$, 
acting on $L^2(M)$ (the book of Davies \cite{14} is a good 
reference).  This 
semigroup maps bounded functions in $L^2(M)$ to bounded functions in 
$L^2(M)$, and the latter determines a semigroup of normal linear 
maps on the abelian von Neumann algebra 
$L^\infty(M)$ which carries nonnegative functions 
to nonnegative functions and fixes the constant functions.   

In order to discuss the symbol of $\Delta$ we introduce local 
coordinates in some open set $U\subseteq M$ to 
identify $U$ with an open region in $\Bbb R^n$.  For clarity, 
we will be explicit with notation.  
At each point $x\in U$ the tangent space $T_x(M)$ is identified 
with $\Bbb R^n$, and for a smooth function $f$ on $M$ the differential 
$df$ takes the following form 
$$
df(x,v)=\frac{d}{dt}f(x+tv)\vert_{t=0}=
\sum_{k=1}^n\frac{\partial f}{\partial x_k}(x)v_k.  
$$
The metric gives rise to a an operator 
function $x\in U\mapsto G(x)$ on $\Bbb R^n$ by way of 
$$
\<v,w\>_{T_x(M)}=\<G(x)v,w\>_{\Bbb R^n}, \qquad v,w\in T_x(M), \qquad x\in U,
$$
where $\<\cdot,\cdot\>_{\Bbb R^n}$ 
denotes the Euclidean inner product on $\Bbb R^n$.  
$G(x)$ is an invertible positive operator on $\Bbb R^n$ 
for every $x\in U$.  
For two vector fields 
$\xi,\eta$ on $M$ we have 
$$
\<\xi(x),\eta(x)\>_{T_x(M)}=\<G(x)\xi(x),\eta(x)\>_{\Bbb R^n}=
\sum_{i,j=1}^n g_{ij}(x)\xi_j(x)\eta_i(x),
$$
for $x\in U$, $(g_{ij}(x))$ being the matrix of $G(x)$ relative 
to the usual orthonormal basis for $\Bbb R^n$.  

The inner product on the tangent space $T_x(M)$ promotes 
naturally to an inner product on the cotangent space 
$T_x^*(M)$.  Indeed, the Riesz lemma implies that every linear functional 
$\rho$ on $T_x(M)$ is associated with a unique vector $\rho_*\in T_x(M)$
via 
$$
\rho(v)=\<v,\rho_*\>_{T_x(M)}, 
$$
and the inner product in $T^*_x(M)$ is defined by 
$$
\<\rho,\sigma\>_{T_x^*(M)}=\<\rho_*,\sigma_*\>_{T_x(M)}.  
$$
With these conventions one finds that for a smooth function $f$ and 
a point $x\in U$, 
$df(x,\cdot)_*$ becomes the vector in $\Bbb R^n$ 
with components $v_1,\dots,v_n$, 
$$
v_i = \sum_{j=1}^ng^{ij}(x)\frac{\partial f}{\partial x_j}(x), 
$$
$(g^{ij}(x))=(g_{ij}(x))^{-1}$ being the matrix of the inverse 
operator $G(x)^{-1}$.  For points $x\in U$ one has 
$$
\<(df)_*,(dg)_*\>_{T_x(M)}=
\sum_{i,j=1}^ng^{ij}(x)
\frac{\partial f}{\partial x_j}\frac{\partial f}{\partial x_i}.  \tag{3.3}
$$

We first recall that the dualized Riemannian metric (whose 
values are inner products on the cotangent spaces 
$T_x^*(M)$) can be linearized naturally so 
that it becomes a $C^\infty(M)$-linear map of the the module $\Omega^{(2)}(M)$
of {\it symmetric} 2-forms.  More explicitly, let $\Omega^1(M)$ be 
the usual module of $1$-forms and let $\Omega^{(2)}(M)$ 
be the submodule of $\Omega^1(M)\otimes_{C^\infty(M)}\Omega^1(M)$ 
consisting of all elements that are fixed under the action of the 
reflection $R$ defined by 
$R: \omega_1\otimes\omega_2\mapsto \omega_2\otimes\omega_1$.  
For $\omega_1,\omega_2\in\Omega^1(M)$ we write 
$\omega_1\omega_2$ for the symmetrized product 
$$
\omega_1\omega_2 = 
\frac{1}{2}(\omega_1\otimes\omega_2+\omega_2\otimes\omega_1)
\in\Omega^{(2)}(M).  
$$
There is a unique homomorphism of $C^\infty(M)$-modules 
$G^*: \Omega^{(2)}(M)\to C^\infty(M)$ satisfying 
$G^*(df\,dg)(x) = \<df,dg\>_{T_x^*(M)}$ for all $x\in M$, 
and in local coordinates (3.3) implies that $G^*$ has the form 
$$
G^*(df\,dg)(x) = 
\sum_{i,j=1}^ng^{ij}(x)\frac{\partial f}{\partial x_j}
\frac{\partial g}{\partial x_i}, \qquad x\in U.  \tag{3.4}
$$
If one knows $G^*$ as a homomorphism of $C^\infty(M)$-modules then 
one also knows the inner product in each cotangent space 
$T_x^*(M)$, and hence one can recover the original metric 
as an inner product on tangent spaces by duality and the 
Riesz lemma as above.  

We now relate these remarks to the symbol of the Laplacian 
$\Delta$ of $M$.  The symbol of any differential operator 
$L:C^\infty(M)\to C^\infty(M)$ of order at most 2 is 
associated with the bilinear form defined on $C^\infty(M)$ by 
$$
f,g\in C^\infty(M)\mapsto L(fg)-fL(g)-gL(f)+fgL(\bold 1).  
$$
A straightforward argument 
shows that there is a (necessarily unique) homomorphism of 
$C^\infty (M)$-modules $\sigma_L: \Omega^{(2)}(M)\to C^\infty(M)$ 
satisfying 
$$
\sigma_L(df\,dg)=L(fg)-fL(g)-gL(f)+fgL(1).  
$$
In particular, this defines the symbol of any 
second order differential operator on $C^\infty(M)$), as an 
element of $\hom(\Omega^2(M),C^\infty(M))$.  

Restricting attention to the operator $L=\Delta$, one sees that 
for each $f\in C^\infty(M)$ the restriction of 
$\Delta(f)$ to $U$ has the form 
$$
\Delta(f)(x)=\sum_{i,j=1}^n g^{ij}(x)
\frac{\partial^2 f}{\partial x_i\partial x_j} + 
\sum_{k=1}^n u_k(x)\frac{\partial f}{\partial x_k}, \tag{3.5}
$$
where $u_1,\dots,u_n$ are appropriate smooth functions 
(see p. 147 of \cite{14}).  
Using (3.5) one easily computes the symbol of $\Delta$, and 
because of the local formula (3.4) for $G^*$ one obtains 
$\sigma_\Delta=2\cdot G^*$.  From these remarks we 
conclude that {\it the symbol of the Laplacian 
(considered as an element of $\hom(\Omega^{(2)}(M),C^\infty(M))$)
is precisely the Riemannian metric of $M$ in its dualized form}.  
\endremark

Returning now to the case of a general CP semigroup 
$\phi=\{\phi_t: t\geq 0\}$ acting on $\Cal B(H)$, let $\Cal A$ be 
the domain algebra of the generator of $\phi$.  Letting $L$  
be the restriction of the generator to $\Cal A$, it is natural to 
ask the extent to which the generator can be identified with something 
analogous to a Riemannian metric (more precisely, to the 
homomorphism of $C^\infty(M)$-modules 
$G^*: \Omega^{(2)}(M)\to C^\infty(M)$ 
that the dualized Riemannian metric determines).  
We have already defined the symbol 
$\sigma_L:\Omega^2(\Cal A)\to\Cal B(H)$ as a homomorphism of 
$\Cal A$-modules, and the following asserts that $\sigma_L$ 
does behave as if it were a (perhaps degenerate) Riemannian metric.  

\proclaim{Proposition}
Let $\phi=\{\phi_t: t\geq 0\}$ be a $CP$ semigroup acting on 
$\Cal B(H)$ and consider the restriction $L$ of the generator 
to the domain algebra $L:\Cal A\to\Cal B(H)$.  Then the symbol 
of $L$ satisfies 
$$
\sigma_L(\omega^*\omega)\leq 0, \qquad \omega\in\Omega^1(\Cal A); 
$$
and more generally for all $\xi_1,\dots,\xi_n\in H$ and 
$\omega_1,\dots,\omega_n\in\Omega^1(\Cal A)$ we have 
$$
\sum_{i,j=1}^n\<\sigma_L(\omega_j^*\omega_i)\xi_i,\xi_j\>\leq 0.  
$$
\endproclaim

The proof is a computation, facilitated by the following formula.  

\proclaim{Lemma}
Let $\omega_1$, $\omega_2$ be elements of $\Omega^1(\Cal A)$ 
having the form 
$$
\omega_k=\sum_{p=1}^s A_{kp}\otimes B_{kp}, \qquad k=1,2,
$$
where $A_{k1}B_{k1}+\dots+A_{ks}B_{ks}=0$ for $k=1,2$.  Then
$$
\sigma_L(\omega_1\omega_2)=-\sum_{p,q=1}^sA_{1p}L(B_{1p}A_{2q})B_{2q}.  
$$
\endproclaim

\demo{proof of Lemma}
Since $\Omega^1(\Cal A)$ is spanned by elements of the form 
$A\cdot dX$, as well as by elements of the form $dY\cdot B$, 
$A,B, X,Y\in\Cal A$, and since $\sigma_L$ is a bimodule homomorphism, 
it suffices to check the formula for $\omega_1$, $\omega_2$ of the 
particular form $\omega_1=dX$, $\omega_2=dY$.  

Writing 
$$
\align
dX\,dY &=(X\otimes\bold 1-\bold 1\otimes X)(Y\otimes \bold 1-\bold 1\otimes Y)\\
&=(X\otimes\bold 1)(Y\otimes\bold 1)-(\bold 1\otimes X)(Y\otimes\bold 1) -
(X\otimes\bold 1)(\bold 1\otimes Y)+(\bold 1\otimes X)(\bold 1\otimes Y) \\
\endalign
$$
the right side of the asserted formula for $\sigma_L(dX\,dY)$ has 
the form
$$
\align
-&(XL(Y)-L(XY)-XL(\bold 1)Y+L(X)Y)=\\
&L(XY)-XL(Y)-L(X)Y+XL(\bold 1)Y=\sigma_L(dX\,dY),
\endalign
$$
as required.\qedd
\enddemo

\demo{proof of Proposition}  
Because of the exact sequence (3.1), every element $\omega\in\Omega^1(\Cal A)$
can be written 
$$
\omega = A_1\otimes B_1+\dots+A_s\otimes B_s, 
$$
where $A_k, B_k$ are elements of $\Cal A$ satisfying 
$A_1B_1+\dots+A_sB_s=0$.  Choose elements 
$\omega_1,\dots,\omega_n\in\Omega^1(\Cal A)$ of the form 
$$
\omega_k=\sum_{p=1}^s A_{kp}\otimes B_{kp}, \qquad k=1,\dots,n
$$
where $\sum_p A_{kp}B_{kp}=0$ for $k=1,\dots,n$.  We have 
$$
\omega_k^*=\sum_{p=1}^nB_{kp}^*\otimes A_{kp}^*
$$
so that 
the product $\omega_k^*\omega_j\in \Omega^2(\Cal A)$ 
is given by 
$$
\omega_k^*\omega_j=
\sum_{p,q=1}^s (B_{kp}^*\otimes A_{kp}^*)(A_{jq}\otimes B_{jq}).  
$$
The Lemma implies that 
$$
\sigma_L(\omega_k^*\omega_j)=
-\sum_{p,q=1}^s B_{kp}^*L(A_{kp}^*A_{jq})B_{jq}. \tag{3.6}
$$

If we now choose vectors $\xi_k\in H$, $k=1,\dots,n$ then we find 
that 
$$
\align
\sum_{k,j}\<\sigma_L(\omega_k^*\omega_j)\xi_j, \xi_k\> &= 
-\sum_{i,j,p,q}\<L(A_{kp}^*A_{jq})B_{jq}\xi_j,B_{kp}\xi_k\>\\
&=-\sum_{\alpha,\beta}\<L(A_\beta^*A_\alpha)\eta_\alpha,\eta_\beta\>, \tag{3.7}
\endalign
$$
where in the third term, $\alpha$ and $\beta$ run over all pairs 
$\{(k,p): 1\leq k\leq n, 1\leq p\leq s\}$ and where the $\eta_\alpha$ are 
defined by 
$\eta_{(k,p)}=B_{kp}\xi_k$.  
Finally, the last term on the right of (3.7) can be rewritten in 
terms of the $n\cdot s\times n\cdot s$ operator matrix 
$\tilde A$ having the entries $A_\alpha$ along a single row and 
zeros along all the other rows as follows 
$$
-\<L^{(n\cdot s)}(\tilde A^*\tilde A)\tilde\eta,\tilde\eta\>, 
$$
where $\tilde\eta$ is the column vector with components $\eta_\alpha$, 
and where $L^{(n\cdot s)}$ is the natural map induced by $L$ on matrices over 
$\Cal A$ be applying $L$ to the elements of the matrix 
term-by-term.  

Thus we have to show that 
$\<L^{(n\cdot s)}(\tilde A^*\tilde A)\tilde\eta,\tilde\eta\>\geq0$.  Recalling that 
the definition of $L$ on elements of $\Cal A$ is 
$$
L(X) = \lim_{t\to 0+}\frac{1}{t}(\phi_t(X)-X)
$$
and the fact that $\Cal A$ is a $*$-subalgebra of the domain of $L$, 
it follows that 
$$
\<L^{(n\cdot s)}(\tilde A^*\tilde A)\tilde\eta,\tilde\eta\>=
\lim_{t\to0+}\frac{1}{t}(\<\phi_t^{(n\cdot s)}(\tilde A^*\tilde A)\tilde\eta,\tilde\eta\>
-\<\tilde A^*\tilde A\tilde\eta,\tilde\eta\>).  
$$

Notice that $\<\tilde A^*\tilde A\tilde\eta,\tilde\eta\>=0$.  Indeed, by 
inspection of the components of the column vector $\tilde A\tilde\eta$ we 
find that it is the column vector having a single (possibly) nonzero component 
and that component is 
$$
\sum_\alpha A_\alpha\eta_\alpha=\sum_{k,p}A_{kp}B_{kp}\xi_k=
\sum_{k=1}^n(\sum_{p=1}^s A_{kp}B_{kp})\xi_k=0,
$$
since $\sum_p A_{kp}B_{kp}=0$ for every $k$.  
Thus we have to show that 
$$
\lim_{t\to0+}\frac{1}{t}
\<\phi_t^{(n\cdot s)}(\tilde A^*\tilde A)\tilde\eta,\tilde\eta\>\geq 0.  \tag{3.8}
$$

Now since for each $t>0$ the map $\phi_t$ is unital and completely 
positive, the Schwarz inequality for completely positive 
maps implies 
$$
\phi_t^{(n\cdot s)}(\tilde A^*\tilde A)\geq 
\phi_t^{(n\cdot s)}(\tilde A)^*\phi_t^{(n\cdot s)}(\tilde A), 
$$
and hence for positive $t$ we have 
$$
\frac{1}{t}\<\phi_t^{(n\cdot s)}(\tilde A^*\tilde A)\tilde\eta,\eta\>\geq 
\frac{1}{t}\<\phi_t^{(n\cdot s)}(\tilde A)\eta,\phi_t^{(n\cdot s)}(\tilde A)\eta\>=
\frac{1}{t}\|\phi_t^{(n\cdot s)}(\tilde A)\eta\|^2.  
$$
We claim that the term on the right tends to zero as $t\to0+$.  
Indeed, since the operator matrix $\tilde A$ belongs to the domain of the 
generator of the CP semigroup 
$\phi^{(n\cdot s)}=\{\phi_t^{(n\cdot s)}: t\geq 0\}$, Lemma 1 implies 
that there is a constant $M>0$ such that for every positive $t$, 
$\|\phi_t^{(n\cdot s)}(\tilde A)-\tilde A\|\leq M\cdot t$.  It follows that 
$$
\|\phi_t^{(n\cdot s)}(\tilde A)\tilde\eta\|=
\|\phi_t^{(n\cdot s)}(\tilde A)\tilde\eta-\tilde A\eta\|\leq 
M\cdot t\cdot\|\tilde\eta\|
$$
and hence 
$$
\limsup_{t\to0+}\frac{1}{t}\|\phi_t^{(n\cdot s)}(\tilde A)\tilde\eta\|^2 
\leq \lim_{t\to0+}\frac{1}{t}(M^2\cdot t^2\cdot\|\tilde\eta\|^2)=0.  
$$
It follows that 
$$
\lim_{t\to0+}\frac{1}{t}
\<\phi_t^{(n\cdot s)}(\tilde A^*\tilde A)\tilde\eta,\tilde\eta\>\geq 
\lim_{t\to0+}\frac{1}{t}\|\phi_t^{(n\cdot s)}(\tilde A)\eta\|^2=0, 
$$
and the inequality (3.8) follows. \qedd
\enddemo

\end